\documentclass[11pt]{amsart}
\usepackage{amsfonts}

\newtheorem{theorem}{\bf Theorem}
\newtheorem{remark}{\bf Remark}[section]
\newtheorem{proposition}{Proposition}[section]
\newtheorem{lemma}{Lemma}

\newtheorem{definition}{Definition}

\usepackage{enumerate,mathrsfs}
\usepackage{amssymb,amsmath,graphicx,amsthm,mathtools,hyperref}

\usepackage{booktabs}

\theoremstyle{plain}
\usepackage{siunitx}

\usepackage{algorithm}
\usepackage{algorithmic,colortbl}

\makeindex
\usepackage[intoc]{nomencl} 

\makenomenclature

\usepackage{hyperref}
\hypersetup{nesting=true,debug=true,naturalnames=true}
\usepackage{graphicx,amssymb,upref}

\usepackage{subcaption}
\usepackage{enumerate,float}
\usepackage{fullpage}

\begin{document}

\title{Hawkes process with tempered Mittag-Leffler kernel}
\author[]{Neha Gupta}
\address{\emph{Department of Mathematics and Statistics,
		Indian Institute of Technology Kanpur, Kanpur 208016, India.}}
\email{nehagpt@iitk.ac.in}
\author[]{Aditya Maheshwari}
\address{\emph{Operations Management and Quantitative Techniques Area, Indian Institute of Management Indore, Indore 453556, India.}}
\email{adityam@iimidr.ac.in}

\begin{abstract}
 In this paper, we propose an extension  of the Hawkes process by incorporating a kernel based on the tempered Mittag-Leffler distribution. This is the generalization of the work presented in \cite{habyarimana2023fractional}.
We derive analytical results for the expectation of the conditional intensity and the expected number of events in the counting process. Additionally, we investigate the limiting behavior of the expectation of the conditional intensity. Finally, we present an empirical comparison of the studied process with its limiting special cases. 
\end{abstract}
\maketitle

\section{Introduction}

Poisson counting processes with time-varying intensities have been found very useful in several applications (see \cite{nhpp-appl1,nhpp-appl2,nhpp-appl3}). If we assume that the arrivals of the counting process are not only dependent on time but also on the (recent) past events, that is, the occurrence of future events for a certain time will be more probable if the last event occurred recently. 
This process was first introduced in the early seventies by A. G. Hawkes (see \cite{hawkes1971spectra, hawkes1971point}) and found applications in several fields, for example, in modelling terrorist activities \cite{Porter_2012}, finance (see 
\cite{Bacry_2015,hainaut2017clustered,Hainaut_2019,hautsch2006modelling, bowsher2007modelling,cartea2014buy,chavez2005estimating}), and
seismology (see \cite{hawkes1973cluster,ogata1988statistical,ogata1998space}).
The Hawkes process (HP)  is a counting process with conditional stochastic intensity  $\{\Lambda(t|\mathcal{H}_t)\}_{t\geq 0}$, where $\mathcal{H}_t$ represents the history of the counting process and is given by 
\begin{equation}
    \label{condint}
\Lambda(t|\mathcal{H}_t) = \Lambda_0 + \alpha \int_0^t f(t-u) \, dN(u),
\end{equation}
\noindent where $\Lambda_0>0$ is the baseline intensity, $\alpha>0$ is the jump size, $f(t)$ is the kernel density function of a positive random variable and $N(t)$ is the counting process. 
In \cite{habyarimana2023fractional}, a Hawkes process of fractional type was introduced by taking the kernel  $f(t)$ as Mittag-Leffler (ML) distribution. Inspired by the formulation in \cite{habyarimana2023fractional}, in this paper, we use the tempered Mittag-Leffler (TML) distribution in place of the ML distribution to define and study the Hawkes process. \\\\
The introduction of TML distribution is motivated by two main reasons. First, the TML generalizes  both the ML and the exponential function and, therefore, offers more flexibility in modelling complex situations (see \cite{temfraccal}). Second, the exponential tempering of the ML function allows us to model both  the long and the short memory processes.  
In  other words, the tempering parameter  can be used to control the memory effects of  past events, which provides greater flexibility for applications. Recently, the fractional Hawkes process (\cite{habyarimana2023fractional}) found application in earthquake shock modelling (see \cite{fhp-earthquake}). In the literature, there is another version of the fractional Hawkes process is known, which is defined as by time-changing the conditional intensity function using an inverse stable subordinator (see \cite{hainaut2020fractional}), its tempered analogue was worked out in \cite{tfhp-gpt-am}. \\

In this paper, we  define the Hawkes process with TML (HPTML) kernel density function and work out the expected intensity, spectral properties and its limiting behaviour. We compare the analytical solution and numerical solution for the expected intensity function. We also empirically compare the HPTML with various processes, eg the Poisson process, fractional Hawkes process and the Hawkes process with the exponential kernel. \\\\
The paper is organized as follows. In Section \ref{sec:prelims}, we present some preliminaries that will be used in this paper. Section \ref{sec:EI} discuss the spectral properties, expected intensity, analytical and numerical solution of the expected conditional intensity function of the HPTML. Finally, we compare various processes with the HPTML in Section \ref{sec:comp}.

\section{Preliminaries}\label{sec:prelims}
In this section, we introduce some definitions, notations and results that will be used later.

 The generalized Mittag-Leffler function, 
 $M_{a, b}^{c}(-\omega y^a)$, introduced by \cite{Mittag-Leffler-general}, is defined as
  \begin{equation}\label{generalML}
       M_{a, b}^{c}(z)=  \sum_{n=0}^{\infty} \frac{(c)_n}{\Gamma(an+b)}\frac{z^n}{n!}, z\in \mathbb{C}
  \end{equation}
  where $a,b,c \in \mathbb{C}$ with real parts of $a,b$ and $c$ are positive and $(c)_n$ is Pochhammer symbol. 
Recall that (see \cite{Mittag-Leffler-general})
\begin{equation}\label{mlintegral}
    \int_{0}^{t}y^{\nu-1}M_{\rho, \mu}^{\nu}(w y^{\rho})(t-y)^{\nu-1}dy= \Gamma(\nu)t^{\nu+\mu-1}M_{\rho, \mu+\nu}^{\nu}(w y^{\rho}).
\end{equation}
The following identity holds true for the two-parameter Mittag-Leffler function \cite[eq. (2.2.1)]{mathai2008special}
\begin{equation} \label{MLidentity}
    M_{a,b}^1 (z) = z M_{a,a+b}^1 (z) + \frac{1}{\Gamma(b)}.
\end{equation}
\\
 A tempered Mittag-Leffler distribution can be introduced by exponential tempering of the Mittag-Leffler distribution (see \cite[eq. (3.24)]{Tempered_Mittag_dis}). Consider $ M_{1,\beta,1-\nu^\beta, \nu}
 (1)$  an exponentially tempered Mittag-Leffler random variable with tempering parameter $\nu \geq 0$. The Laplace transform $\mathcal{L}(\cdot)$ of the probability density function $f_{\beta, \nu}$ is given by (see \cite[eq. (3.25)]{Tempered_Mittag_dis})
 \begin{equation}
\label{ltmlrv}
\mathcal{L}[f_{\beta, \nu}(t): s]= \tilde{f}_{\beta, \nu} (s) = \int_0^\infty \mathrm{e}^{-st} f_{\beta, \nu} (t) \; dt = \frac{1}{1-\nu^\beta +(\nu+s)^\beta},
\end{equation}
where $\nu>0$, $\beta \in (0, 1)$ and $\mathbb{R}(s)>0$.
\begin{lemma}[Hardy and Littlewood 1930]
\label{tauberian1}
If $f(t)$ is positive and integrable over every finite range $(0,T)$ and $\mathrm{e}^{-st} f(t)$ is integrable over $(0, \infty)$ for every $s >0$ and if
\begin{equation*}
\tilde{f} (s) = \int_0^\infty \mathrm{e}^{-st} f(t) \, dt  \sim H s^{-a}
\end{equation*}
for $s \to 0$ with $a >0$, $H>0$, then as $t \to \infty$
\begin{equation*}
F(t) = \int_0^t f(u) \, du \sim \frac{H}{\Gamma(1+a)} t^a.
\end{equation*}
\end{lemma}
\begin{definition}[Hawkes Process]
    Consider $(N(t) : t \geq 0)$, (see \cite{ hawkes1971point, hawkes1971spectra}) a counting process, with associated history
$\{\mathcal{H}_t\}_{t \geq 0}=  (\mathcal{H}(t) : t \geq 0)$, that satisfies the following 
$$
\mathbb{P}(N(t+h)-N(t)=n|\mathcal{H}_t)
     \begin{cases}
        \Lambda(t|\mathcal{H}_t) h+o(h), & n= 1,\\
                  1- \Lambda(t|\mathcal{H}_t) h +o(h) , & n=0,\\
                  o(h), & n\geq 2.
           \end{cases}
$$
The conditional intensity function is defined in \eqref{condint}.
\end{definition}
The trivial case $\alpha=0$, is a homogeneous Poisson process with intensity rate $ \Lambda_0$.

\section{Expected intensity and limiting behavior}
\label{sec:EI}

In this section, we  work out analytical results related to the expected intensity of the Hawkes process with TML intensity kernel. We compute the Laplace transform (LT) of the expected intensity and then invert the LT to get the desired result. Moreover, we  find out the asymptotic behaviour of the expected intensity. The mean of the number of events is also discussed.

\noindent Let $\lambda (t)$ be the expected intensity of the conditional intensity function $\Lambda(t|\mathcal{H}_t)$ is given by
\begin{equation*}
\lambda (t) = \mathbb{E} [\Lambda(t|\mathcal{H}_t)].
\end{equation*}
Using \eqref{condint}, we have that 
\begin{equation*}
\lambda(t) = \Lambda_0 +\alpha \int_0^t f_{\beta, \nu} (t -u) \lambda(u) \; du.
\end{equation*}
Now, using \eqref{ltmlrv}, we get the LT of this equation
\begin{equation*}
\tilde{\lambda}(s)= \frac{\Lambda_0}{s
(1-\alpha \tilde{f}_{\beta, \nu}(s))}=\frac{\Lambda_0}{s}\left[\frac{1-\nu^\beta+(\nu+s)^\beta}{1-\alpha-\nu^\beta+(\nu+s)^\beta}\right] .
\end{equation*}
Changing the time scale using the factor $\gamma > 0$, the above equation reduces to
\begin{equation}
\label{expectedintensity}		  
\tilde{\lambda}(s)=\frac{\Lambda_0}{s}\left[\frac{\gamma-\nu^\beta+(\nu+s)^\beta}{\gamma(1-\alpha)-\nu^\beta+(\nu+s)^\beta}\right].
\end{equation}

\begin{figure}[!ht]
			\centering
			\begin{minipage}{.5\textwidth}
				\centering
				\includegraphics[width=1\linewidth]{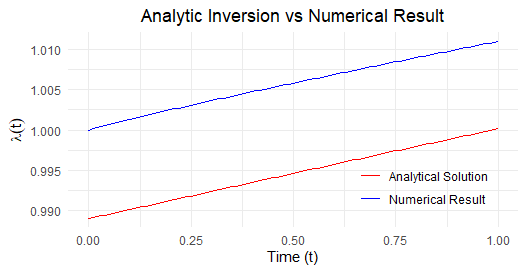}
			\end{minipage}%
			\begin{minipage}{.5\textwidth}
				\centering
				\includegraphics[width=1\linewidth]{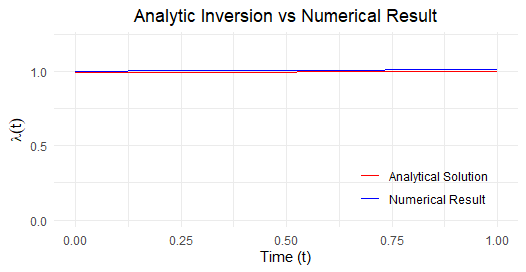}
			\end{minipage}

			\caption{Comparison between the exact formula of expected intensity \eqref{analyticalinverse}
for $\beta = 0.9$ and $\nu=1$ (red line) and the numerical inversion of the Laplace transform in \eqref{expectedintensity} (blue line) for the parameters $\Lambda=1, \alpha=0.1, \gamma=0.1$. \label{fig:comp} }
		\end{figure}

\noindent We now state the following Theorem.
\begin{theorem} Let $\beta \in (0,1)$ and $\nu >0 $, the analytical inverse to equation \eqref{expectedintensity} 
is given by
\begin{align}\label{analyticalinverse}
\lambda(t) &=\Lambda_0\left[\frac{\nu^\beta-\gamma}{(\alpha-1)\gamma+\nu^\beta} +\frac{\alpha \gamma e^{-\nu t}}{(\alpha-1)\gamma+\nu^\beta}\sum_{m=0}^{\infty} \nu^m t^{m} M_{\beta,m+1}^1(((\alpha-1)\gamma+\nu^\beta)t^{\beta})\right],
\end{align}
where $M^1_{\beta,m+1}(z)$ is the two-parameter Mittag-Leffler function defined in \eqref{generalML}.
\end{theorem}
\begin{proof}
To prove the result, we first rewrite \eqref{expectedintensity} as
$$
\tilde{\lambda}(s) = \frac{\Lambda_0}{s} \left[1+\frac{\gamma \alpha}{(1-\alpha)\gamma-\nu^\beta + (s+\nu)^\beta} \right].
$$
 The inverse LT of the above equation,
 \begin{align}\label{InLT}
     \mathcal{L}^{-1}[\tilde{\lambda}(s): t] &= \Lambda_0  \left[\mathcal{L}^{-1}\left[\frac{1}{s} :t\right] +\mathcal{L}^{-1}\left[\frac{\gamma\alpha}{s((1-\alpha)\gamma-\nu^\beta + (s+\nu)^\beta)} :t\right]\right].
 \end{align}
 We first calculate the second term of the above equation. Let us denote it by $G(s)$, which is equal to
$$
G(s)=\frac{\alpha\gamma}{(1-\alpha)\gamma-\nu^\beta + (s+\nu)^\beta}.
$$
Using the shifting property  $\mathcal{L}[e^{at}f(t): s]=\tilde{f}(s-a)$, where $\tilde{f}(s)=\mathcal{L}[f(t): s]$,  the expression simplifies to
$$
\mathcal{L}^{-1}[G(s): t]= e^{-\nu t}\mathcal{L}^{-1}\left[\frac{\gamma\alpha}{(1-\alpha)\gamma-\nu^\beta + s^\beta}\right].
$$
We know that (see \cite{Mittag-Leffler-general})
$$
 \mathcal{L}^{-1}\left[\frac{s^{ac-b}}{(s^a+\omega)^c}: t\right]= t^{b-1}M_{a, b}^{c}(-\omega t^a).
 $$
As a result, the inverse LT of $G(s)$ is
 $$
 \mathcal{L}^{-1}[G(s): t]=\alpha\gamma  e^{-\nu t} t^{\beta-1}M_{\beta, \beta}^1(((\alpha-1)\gamma+\nu^\beta )t^{\beta}).
 $$
 Moreover,
 \begin{align}\label{Gs}
\mathcal{L}^{-1}\left[\frac{G(s)}{s}: t\right]
  &= \alpha\gamma\int_{0}^t e^{-\nu u} u^{\beta-1}M_{\beta, \beta}^1(((\alpha-1)\gamma+\nu^\beta)u^{\beta}) du\nonumber \\
  &=\alpha\gamma e^{-\nu t}\sum_{m=0}^{\infty}\nu^m t^{\beta+m}M_{\beta, \beta+m+1}^1(((\alpha-1)\gamma+\nu^\beta)t^{\beta}) \text{ (using \eqref{mlintegral}}.
  \end{align}
Substituting \eqref{Gs} in \eqref{InLT}, we get
\begin{equation*}
\mathcal{L}^{-1}[\tilde{\lambda}(s): t] = \Lambda_0\left[1+\alpha\gamma e^{-\nu t}\sum_{m=0}^{\infty}\nu^m t^{\beta+m}M_{\beta, \beta+m+1}^1(((\alpha-1)\gamma+\nu^\beta)t^{\beta})\right]
\end{equation*}
Using \eqref{MLidentity} with $a = \beta$, $b=m+1$ and $z=((\alpha-1)\gamma+\nu^\beta)t^{\beta})$, we obtain the following expression.
\begin{align*}
\lambda(t) = \Lambda_0 \left[1-\frac{\gamma \alpha}{(\alpha-1)\gamma+\nu^\beta}+\frac{\alpha\gamma e^{-\nu t}}{(\alpha-1)\gamma+\nu^\beta}\sum_{m=0}^{\infty}
\nu^m t^{m} M_{\beta, m+1}^1(((\alpha-1)\gamma+\nu^\beta)t^{\beta})\right].
\end{align*}
Rearrange the above expression to get the desired result in \eqref{analyticalinverse}.
\end{proof}
\noindent We compare the results for the analytical inverse with the numerical inverse of the LT of the expected intensity in Figure \ref{fig:comp}. In the left pane of the figure, the difference is highlighted, and it is of the order of ~0.01. 

\begin{remark} The process considered above generalizes of the fractional Hawkes process (see \cite{habyarimana2023fractional}).  
    Let $\nu =0$, then inverse LT  reduces to 
    $$
    \tilde{\lambda}(s)=\frac{\Lambda_0}{s}\left[\frac{\gamma+s^\beta}{\gamma(1-\alpha)+s^\beta}\right],
    $$
    which is equal to the LT of the expected intensity of the fractional Hawkes process defined in \cite{habyarimana2023fractional}.
\end{remark}

\noindent We now discuss the asymptotic results for the LT of the expected intensity in the following Lemma.
\begin{lemma}
\label{asymptotic5}
Consider $\tilde{\lambda} (s)$ given by equation \eqref{expectedintensity}. As $s\to0$, we have that
\begin{equation*}
\tilde{\lambda} (s) \sim \frac{\Lambda_0(\gamma-\nu^\beta)s^{-1}}{\gamma(1-\alpha)-\nu^\beta}\left[\sum_{n=0}^\infty \left(\frac{\nu^{\beta n}}{((\alpha-1)\gamma+\nu^\beta)^n}+\frac{\nu^{\beta(n+1)}}{((\alpha-1)\gamma+\nu^\beta)^n (\gamma-\nu^\beta)^n}\right)\right].
\end{equation*}
\end{lemma}
\begin{proof}
We refer to the Puiseux-Newton series (see \cite{puiseux1850recherches}) to expand eq \eqref{expectedintensity} around $s =0$ to get the following expression:
\begin{align*}
\tilde{\lambda} (s) &=\frac{\Lambda_0(\gamma-\nu^\beta)s^{-1}}{\gamma(1-\alpha)-\nu^\beta}\left[\left(1+\frac{(\nu+s)^\beta}{\gamma-\nu^\beta}\right)\left(1+\frac{(s+\nu)^\beta}{\gamma(1-\alpha)-\nu^\beta}\right)^{-1}\right]\\
&=\frac{\Lambda_0(\gamma-\nu^\beta)s^{-1}}{\gamma(1-\alpha)-\nu^\beta}\left[\left(1+\frac{(\nu+s)^\beta}{\gamma-\nu^\beta}\right)\left(\sum_{n=0}^{\infty}\frac{(\nu+s)^\beta}{(\gamma(\alpha-1)+\nu^\beta)^n}\right)\right]\\
&= \frac{\Lambda_0(\gamma-\nu^\beta)s^{-1}}{\gamma(1-\alpha)-\nu^\beta}\left[\sum_{n=0}^\infty \left(\frac{\nu^{\beta n}}{((\alpha-1)\gamma+\nu^\beta)^n}+\frac{\nu^{\beta(n+1)}}{((\alpha-1)\gamma+\nu^\beta)^n (\gamma-\nu^\beta)^n}\right)\right]\nonumber\\
&+ \sum^{\infty}_{k=1}\frac{\Lambda_0(\gamma-\nu^\beta)s^{k-1}}{\gamma(1-\alpha)-\nu^\beta}\left[\sum_{n=0}^\infty\left({n\beta \choose k}\frac{\nu^{\beta n-k}}{((\alpha-1)\gamma+\nu^\beta)^n}+{(n+1)\beta \choose k}\frac{\nu^{\beta(n+1)-k}}{((\alpha-1)\gamma+\nu^\beta)^n (\gamma-\nu^\beta)^n}\right)\right].
\end{align*}
This implies that
\begin{equation*}
\lim_{s \to 0} \frac{\tilde{\lambda}(s)}{\frac{\Lambda_0(\gamma-\nu^\beta)s^{-1}}{\gamma(1-\alpha)-\nu^\beta}\left[\sum_{n=0}^\infty \left(\frac{\nu^{\beta n}}{((\alpha-1)\gamma+\nu^\beta)^n}+\frac{\nu^{\beta(n+1)}}{((\alpha-1)\gamma+\nu^\beta)^n (\gamma-\nu^\beta)^n}\right)\right]} = 1. \qedhere
\end{equation*}
\end{proof}

\noindent Using Lemma \ref{asymptotic5}, we have the following Theorem.
\begin{theorem}
Let $\lambda(t) = \mathcal{L}^{-1} \{ \tilde{\lambda} (s) \} (t)$ with $ \tilde{\lambda} (s)$ given by equation \eqref{expectedintensity}. Then $$\lim_{t \to \infty} \lambda(t) = \frac{\Lambda_0(\gamma-\nu^\beta)}{\gamma(1-\alpha)-\nu^\beta}\left[\sum_{n=0}^\infty \left(\frac{\nu^{\beta n}}{((\alpha-1)\gamma+\nu^\beta)^n}+\frac{\nu^{\beta(n+1)}}{((\alpha-1)\gamma+\nu^\beta)^n (\gamma-\nu^\beta)^n}\right)\right].$$
\end{theorem}
\begin{proof}
Let $H = \frac{\Lambda_0(\gamma-\nu^\beta)}{\gamma(1-\alpha)-\nu^\beta}\left[\sum_{n=0}^\infty \left(\frac{\nu^{\beta n}}{((\alpha-1)\gamma+\nu^\beta)^n}+\frac{\nu^{\beta(n+1)}}{((\alpha-1)\gamma+\nu^\beta)^n (\gamma-\nu^\beta)^n}\right)\right]$ and
$a = 1$, we the expected number of events asymptotically behaves in the following manner (using Lemma \ref{tauberian1} and Lemma \ref{asymptotic5}) for large $t$
\begin{equation}
\label{expnasympt}
\mathbb{E} [N(t)] = \int_0^t \lambda (u) \,du \sim \frac{\Lambda_0(\gamma-\nu^\beta)}{\gamma(1-\alpha)-\nu^\beta}\left[\sum_{n=0}^\infty \left(\frac{\nu^{\beta n}}{((\alpha-1)\gamma+\nu^\beta)^n}+\frac{\nu^{\beta(n+1)}}{((\alpha-1)\gamma+\nu^\beta)^n (\gamma-\nu^\beta)^n}\right)\right] t.
\end{equation}
Differentiating  \eqref{expnasympt} leads to 
\begin{equation*}
\lambda (t) \sim \frac{\Lambda_0(\gamma-\nu^\beta)}{\gamma(1-\alpha)-\nu^\beta}\left[\sum_{n=0}^\infty \left(\frac{\nu^{\beta n}}{((\alpha-1)\gamma+\nu^\beta)^n}+\frac{\nu^{\beta(n+1)}}{((\alpha-1)\gamma+\nu^\beta)^n (\gamma-\nu^\beta)^n}\right)\right],
\end{equation*}
which completes the proof. \qedhere\end{proof}
\begin{remark}
For $\nu=0$, the asymptotic behavior of the expected intensity for HPTML is the same fractional HP for large $t \to \infty$ (see \cite{habyarimana2023fractional}),
    $$
    \lambda(t) \sim \frac{\Lambda_0}{(1-\alpha)}.
    $$
\end{remark}
\noindent The expected number of events, $\mathbb{E}[N(t)]$ up to time $t$ is given by 
\begin{align*}
\mathbb{E}[N(t)] &= \Lambda_0\left[\frac{\nu^\beta-\gamma}{(\alpha-1)\gamma+\nu^\beta}t +\frac{\alpha \gamma e^{-\nu t}}{(\alpha-1)\gamma+\nu^\beta}\right.\\
&\left.\times\sum_{m=0}^{\infty}\sum_{n=0}^{\infty}\nu^{n+m} t^{m+n+1} M_{\beta,m+n+2}^1(((\alpha-1)\gamma+\nu^\beta)t^{\beta})\right],
\end{align*} 
where
$$
\int_0^{t} e^{-\nu y}y^{m}M_{\beta, m+1}^1(((\alpha-1)\gamma+\nu^\beta)y^{\beta})dy= e^{-\nu t}\sum_{n=0}^{\infty}\nu^n t^{m+n+1}M_{\beta, m+n+2}^1(((\alpha-1)\gamma+\nu^\beta)y^{\beta}).
$$

\section{Comparing HPTML distribution with existing point processes}\label{sec:comp}
In this section, we will use the Monte-Carlo simulation approach to compare the distribution of the number of events of the HPTML with the following processes:
\begin{enumerate}[(a)]
    \item Poisson process, 
    \item Hawkes process with exponential exciting function, and
     \item Fractional Hawkes process (with ML kernel).

\end{enumerate}

\subsection*{Special cases}
The Hawkes process with TML kernel generalizes some of the processes available in the literature. Taking  particular values of parameters such as the tempering parameter $ \nu\geq 0$, fractional index/branching ratio $0<\beta\leq 1$ and jump size $\alpha \geq 0$, we get the following known processes:
\begin{itemize}
    \item Let $\beta=1,\; \alpha=0$ and $\nu=0$, then this counting Hawkes process behaves same as Poisson process. We demonstrate it empirically in Figure \eqref{fig:TFHP-PP-eq}.
    \item Let $\beta=1$ and $\nu=0$. The the conditional intensity function \eqref{condint}  of the Hawkes process reduces to the intensity function of the Hawkes process with exponential exciting function (see Figure \eqref{fig:TFHP-EXP}).
    \item Let $\nu=0$, the  considered process reduces to the fractional Hawkes process of \cite{habyarimana2023fractional} (see Figure \ref{fig:TFHP-FHP-eq}. 
\end{itemize}
\subsection*{Discussion on graphs} 
\begin{figure}[!ht]
    \centering
    \includegraphics[width=0.7\linewidth]{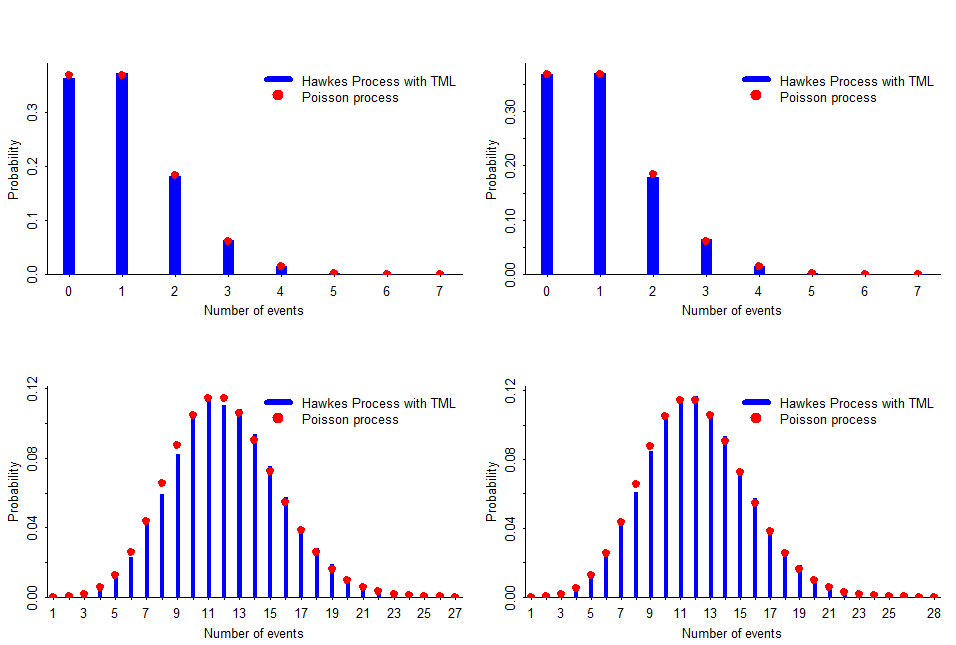}
    \caption{Comparison between the distribution for HPTML (vertical blue bars) and for limiting the Poisson process with the
exponential kernel (red dots). The parameters are $\Lambda_0=1,\; \gamma=1,\; \alpha=0.01,\;  \nu=0.01,\;  t=1, 15, $ from top to bottom and $\beta=0.9,\; 0.99$ from left to right.}
    \label{fig:TFHP-PP-eq}
\end{figure}

\begin{figure}[!ht]
    \centering
    \includegraphics[width=0.7\linewidth]{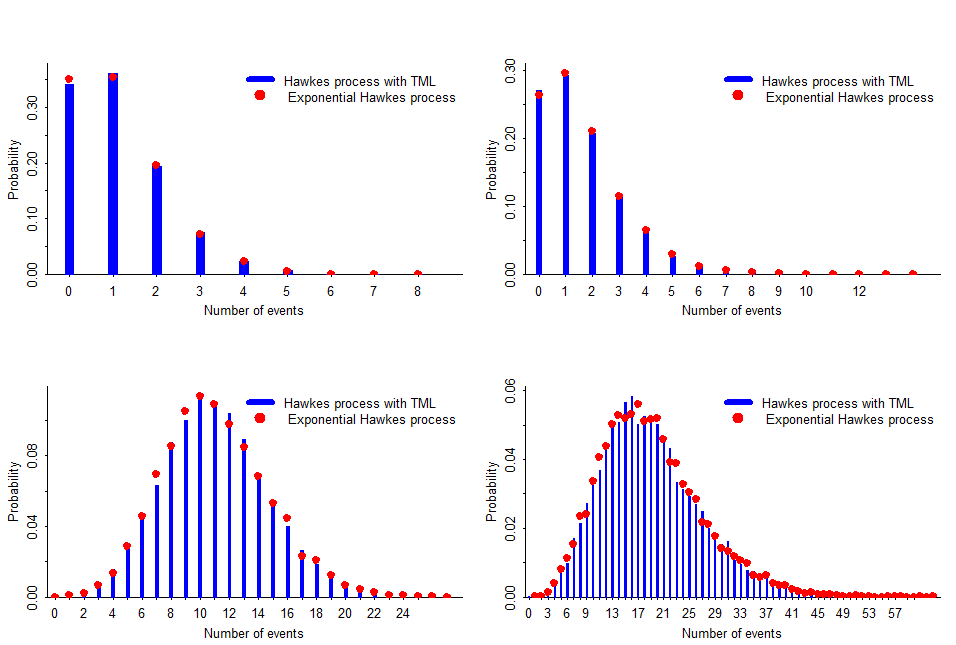}
    \caption{ Comparison between the distribution for HPTML (vertical blue bars) and for limiting the Hawkes process with the
exponential kernel (red dots). The parameters are $\Lambda_0 = 1, \gamma = 1, \beta = 0.99, \nu=0.01,  t = 1, 10$ from top to bottom and $\alpha = 0.1, 0.5$ from left to right.}
    \label{fig:TFHP-EXP}
\end{figure}

\begin{figure}
    \centering
    \includegraphics[width=0.7\linewidth]{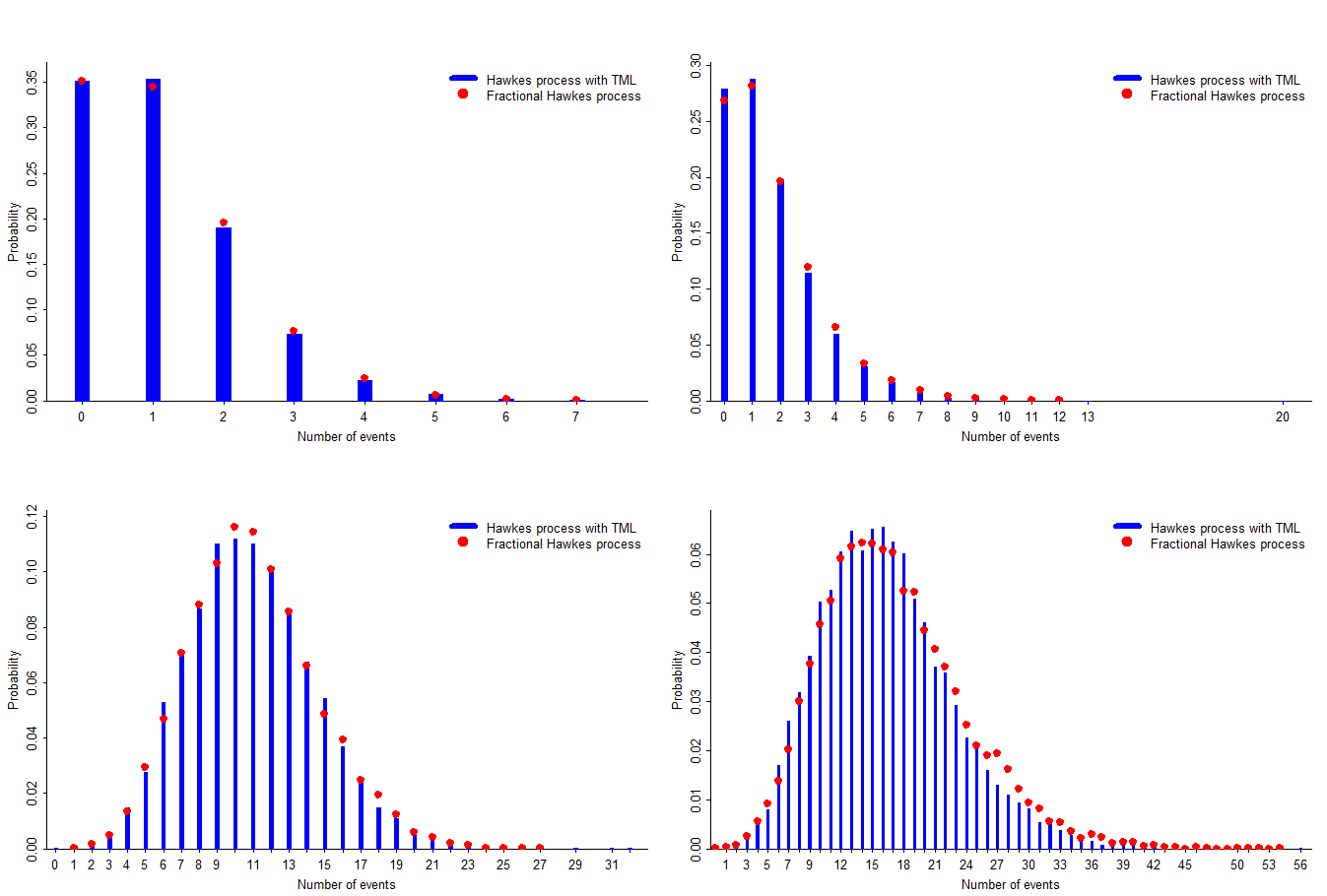}
    \caption{Comparison of the distribution of HPTML (vertical blue bar) and the limiting fractional Hawkes process with the Mittag-Leffler kernel (red dots). The parameters are $\Lambda_0=1, \gamma=1, \beta=0.99, , \nu=0.01, \; t=1,10$ from top to bottom and $\alpha=0.1, \;0.5$ from left to right.}
    \label{fig:TFHP-FHP-eq}
\end{figure}

\begin{figure}[!ht]
    \centering
    \includegraphics[width=0.8\linewidth]{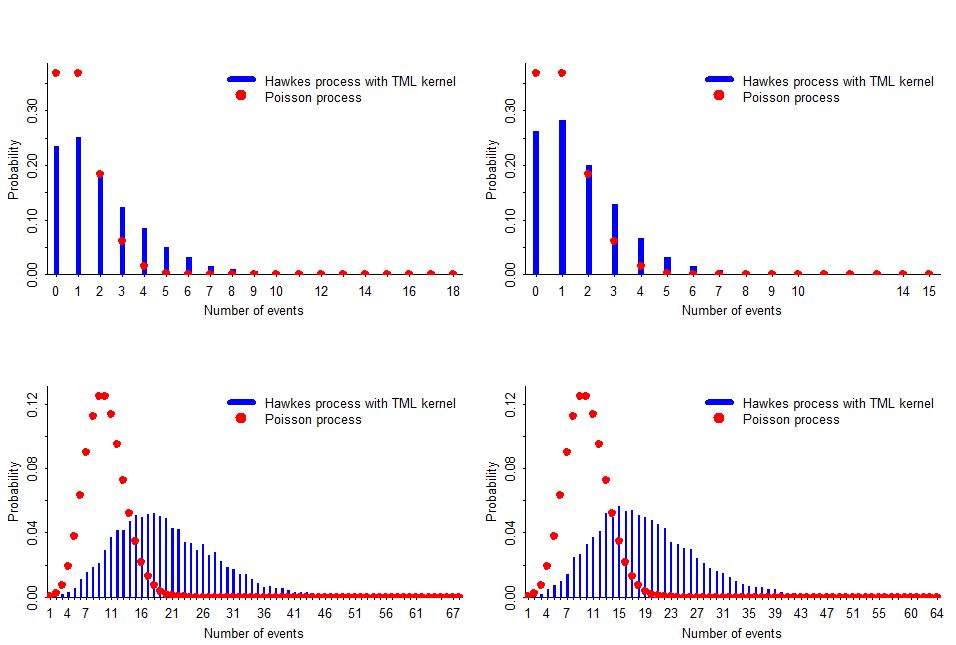}
    \caption{Comparison between the distribution for HPTML (vertical blue bars) and the
Poisson process (red dots). The parameters are $\alpha=0.5,\; \Lambda_0=1, \nu=1, \gamma=1, t=1,10$ from top to bottom and $\beta = 0.5, 0.7$ from left to right.}
    \label{fig:HP-PP-neq}
\end{figure}

We present some graphs to compare the limiting cases of the HPTML. The codes used to make these graphs are inspired from \cite{habyarimana2023fractional}. In Figure \ref{fig:TFHP-FHP-eq}, we plot the HPTML for a small value of $\nu=0.01$ and observe that its distribution matches with the distribution of the fractional Hawkes process of \cite{habyarimana2023fractional}. Figure \ref{fig:TFHP-EXP} shows that when we take $\beta=0.99$ and $\nu=0.01$, the HPTML matches it with its limit case Hawkes process with the exponential kernel. Similarly, in Figure \ref{fig:TFHP-PP-eq}, we take the limiting case when $\nu$ is close to zero,  $\beta$ is close to one and $\alpha$ close to zero to observe the empirical limiting distribution of the HPTML with the Poisson process.  Also, we can observe the deviation of the HPTML with the Poisson process in non-limiting cases in Figure \ref{fig:HP-PP-neq}.\\

\clearpage

\end{document}